\documentclass[10pt]{article}

\usepackage{amsmath,amsthm,amssymb,verbatim}
\usepackage[letterpaper,margin=3cm]{geometry}
\usepackage{graphicx}
\usepackage[ps2pdf,backref,colorlinks]{hyperref}

\usepackage{color}

\newcommand{\la}{\langle}
\newcommand{\ra}{\rangle}

\long\def\comment#1\endcomment{}
\def\G{{\mathcal G}}
\def\Z{{\mathbb Z}}

\newtheorem{theorem}{Theorem}
\newtheorem{proposition}[theorem]{Proposition}
\newtheorem{conjecture}{Conjecture}
\newtheorem{corollary}[theorem]{Corollary}
\newtheorem{lemma}[theorem]{Lemma}

\newenvironment{itemize*}
  {\begin{itemize}
     \setlength{\itemsep}{-3pt}
     \setlength{\parsep}{0pt}}
  {\end{itemize}}

\begin{document}

\title{Critical Percolation of Free Product of Groups}
\author{Iva Koz\'akov\'a\footnote{Electronic version of an article published as 
International Journal of Algebra and Computation,
Volume No.18, Issue No.4, June 2008, Page: 683 - 704,
DOI:10.1142/S0218196708004524
\copyright World Scientific Publishing Company, 
http://ejournals.wspc.com.sg/ijac/ijac.shtml}}
\date{}
\maketitle

\begin{abstract}
In this article we study percolation on the Cayley graph of a free
product of groups.

The critical probability $p_c$ of a free product
$G_1*G_2*\cdots*G_n$ of groups is found as a solution of an equation
involving only the expected subcritical cluster size of factor
groups $G_1,G_2,\dots,G_n$. For finite groups these equations are
polynomial and can be explicitly written down. The expected
subcritical cluster size of the free product is also found in terms
of the subcritical cluster sizes of the factors. In particular, we
prove that $p_c$ for the Cayley graph of the modular group
$\hbox{PSL}_2(\mathbb Z)$ (with the standard generators) is
$.5199...$, the unique root of the polynomial
$2p^5-6p^4+2p^3+4p^2-1$ in the interval $(0,1)$.

In the case when groups $G_i$ can be ``well approximated'' by a
sequence of quotient groups, we show that the critical probabilities
of the free product of these approximations converge to the critical
probability of $G_1*G_2*\cdots*G_n$ and the speed of convergence
is exponential. Thus for residually finite groups, for example, one
can restrict oneself to the case when each free factor is finite.

We show that the critical point, introduced by Schonmann,  $p_{\mathrm{exp}}$ of the free
product is just the minimum of $p_{\mathrm{exp}}$ for the factors.

\comment
We also find explicit formulas for the right derivative of
percolation function and the left derivative of the expected cluster
size function of a free product at $p_c$. It immediately implies the
known result that the values of critical exponents of a free product
are mean-field.
\endcomment
\end{abstract}

%\tableofcontents

\section{Introduction}

\subsection{Percolation}

We will use the notation $\mathcal{G}=(V,E)$ for a graph with the
vertex set $V$ and the edge set $E$. A graph $\mathcal{G}$ is said
to be {\em locally finite} if each vertex has finitely many
neighbors, and {\em transitive} if for any two vertices $u$, $v$ in
$V$ there is an automorphism of $\mathcal{G}$ mapping $u$ to $v$.

An edge of the graph is called a  {\em bond}. A {\em Bernoulli bond
percolation} on $\mathcal{G}$ is a product probability measure
$\mathrm{P}_p$ on the space $\Omega=\{0,1\}^E$, the subsets of the
edge set $E$. For any {\em realization} $\omega \in \Omega$, the
bond $e\in E$ is said {\em open} if $\omega(e)=1$ and {\em closed}
otherwise. For $0\leq p\leq 1$ the product measure is defined via
$\mathrm{P}_p(\omega(e)=1)=p$ for all $e\in E$. Thus each bond is
open with probability $p$ independently of all other bonds. We write
$\mathrm{E}_p$ for the {\em expected value} with respect to
$\mathrm{P}_p$.

For any realization $\omega$, open edges form a random subgraph of
$\mathcal{G}$. An {\em(open) cluster} is a connected component of such
subgraph $\omega$.
An open cluster containing the origin is denoted by $C$ and the number of vertices in $C$ by $|C|$.  {\em Percolation
function} is defined to be the probability that the origin is contained in an infinite cluster, i.e. $\theta(p)=\mathrm{P}_p(|C|=\infty)$.

Depending on the parameter $p$ every subgraph $\omega$ has either no
infinite cluster, or infinitely many infinite clusters ({\em
non-uniqueness phase}), or only one infinite cluster ({\em
uniqueness phase}) $\mathrm{P}_p$-almost surely. H\"aggstr\"om and
Peres \cite{haggstrom:mono} have shown that for transitive graphs
there are two phase transition values of $p$; $p_c$ and $p_u$, such
that for $0\leq p < p_c$ all clusters are finite, non-uniqueness
phase occurs for $p_c< p < p_u$  and if $p_u < p \leq1$ there is
unique infinite cluster $\mathrm{P}_p$-a.s. The critical probability
is then equivalently defined by
$$
p_c=\inf\{p:\theta(p)>0\}.
$$

There is another critical value of $p$ which can be defined based on
the probability of open path between two vertices.
\begin{align*}
p_{\mathrm{exp}}=\sup\{p : \exists_{C,\gamma>0} \forall_{x,y \in V}
\mathrm{P}_p(x\leftrightarrow y)\leq C e^{-\gamma
\mathrm{dist}(x,y)} \}
\end{align*}
This critical points lies between $p_c$ and $p_u$ as pointed out by
Schonmann \cite{schonmann:nonamen}.

As $p$ approaches $p_c$ the behavior of the percolation function and mean cluster size is studied. Assume $\theta(p)$ is continuous at $p_c$ and that
\begin{align*}
\theta(p)&\approx(p-p_c)^\beta\quad \text{ as }p \searrow p_c,\\
\mathrm{E}_p(|C|)&\approx(p_c-p)^\gamma\quad \text{ as } p \nearrow p_c.
\end{align*}
Then we say that $\beta$ and $\gamma$ are {\em critical exponents}.

The {\em Cayley graph} of a group $G$ with respect to the finite set of
generators $S$ is the graph $\mathcal{G}$ with vertices $V=G$ and
$\{g,h\}\in E$ iff $g^{-1}h\in S$ (with the appropriate multiplicity).
This graph is always locally finite and transitive.

Percolation characteristics of a Cayley graph ($p_c, \beta, \gamma$,
etc.) of the group are important invariants of the Cayley graph and
the group related to the spectral radius, $l_2$-Betti numbers, the
Cheeger constant, amenability, etc.

For example if a group is amenable then the non-uniqueness phase is
empty, i.e. $p_c=p_u$. On the other hand Pak and Smirnova-Nagnibeda
\cite{pak:nonuniq} showed that if $G$ is non-amenable then there is
a generating set $S$ of $G$ such that the percolation on a Cayley
graph of $G$ with respect to $S$ has nontrivial non-uniqueness
phase. The general problem whether this is true for all Cayley
graphs of non-amenable groups is still open.

Recall that the {\em Cheeger constant} of a graph $\G(V,E)$
is defined by
$$h(\G)=\inf_{K}\frac{|\partial K|}{|K|},$$
where $K$ is a finite subset of $V$ and $\partial K$, the boundary
of $K$, contains all edges in $E$ with exactly one endpoint in $K$.

There are several general inequalities involving $p_c$. The critical
probability $p_c$ of a quotient graph  does not exceed the $p_c$ of the original graph (see Campanino
\cite{campanino}). In particular, the $p_c$ of a Cayley graph of any
factor group of $G$ is at most the $p_c$ of a Cayley graph of $G$
itself (with respect to the corresponding generating sets).

It is easy to show (using the expected cluster size) that
\begin{align}
    p_c\geq \frac{1}{2|S|-1}, \label{eq:est1}
\end{align}
where the equality holds for free groups
(i.e. when the Cayley graph is a tree).
On the other hand Benjamini and Schramm \cite{benjamini:beyond} proved that
\begin{align}
    p_c\leq \frac{1}{h(\G)+1} \label{eq:est2}
\end{align}
and again the equality holds for free groups.

Gaboriau \cite{gaboriau:harmonic} related  harmonic Dirichlet
functions on a graph to those on the infinite clusters in the
uniqueness phase. He also proved that the first $\ell^2$-Betti
number of a group does not exceed $\frac12 (p_u-p_c)$.

Note also that random subgraphs of the Cayley graph are crucial in
the study of generic properties of a group and the average case
complexity of the word problem \cite{MyasnikovGurevich}.

Probabilistic properties of Cayley graphs (say, properties of the
random walks) have been extensively studied \cite{woess:walks}, but properties of percolation initiated by
Benjamini and Schramm \cite{benjamini:beyond} have not been studied
as much mostly because it is usually difficult to find the explicit
values of the percolation characteristics even for relatively simple
graphs.

Explicit values of $p_c$ are known only for some special cases. For
example, for lattices in $\mathbb{R}^2$ the value of $p_c$ is
obtained using  dual graphs (for $\mathbb{R}^d$ with $d\ge 3$ the
values of $p_c$ are not known). For square lattice, Kesten (\cite{kesten}) proved $p_c=1/2$, for
triangular lattice $p_c=2\sin(\pi/18)$, and for hexagonal lattice
$p_c=1-2\sin(\pi/18)$ (see Grimmett \cite{grimmett:book}). Ziff and
Scullard \cite{ziff:twodim} recently found $p_c$ for a larger class
of lattices in $\mathbb{R}^2$ (they considered graphs which can be
decomposed onto certain self-dual arrangement). Grimmett and Newmann
\cite{grimmett:dimension} studied the percolation on a direct
product of regular tree with $\mathbb{Z}$, they discuss how $p_c$
and $p_u$ changes with the degree of the tree. Lyons studied percolation on arbitrary trees \cite{lyons:book}.

Note that many of the graphs from the previous paragraph are Cayley
graphs of groups: the square lattice in ${\mathbb R}^2$, the
triangular lattice in ${\mathbb R}^2$, some trees, and the direct
product of a tree and ${\mathbb Z}$.

The critical exponents are known for some lattices in
$\mathbb{R}^2$, for $\mathbb{R}^d$, $d\ge 19$, trees and more
generally Cayley graphs with infinitely many ends. For last three
examples $\beta=1, \gamma=-1$.

In this article we will focus on the Cayley graphs of free products
of groups $G_1*G_2*\cdots*G_n$. Some probabilistic properties of
such graphs have been studied before. For example Mairesse and

Math\'eus \cite{mairesse:walk} considered the transient
nearest-neighbor random walk on the Cayley graph of free product of
finite groups. We obtain, among other results, explicit formulas for
$p_c$ (as solutions of some equations), for the Schonmann's critical
point $p_{\mathrm{exp}}$, and for the right derivative of the percolation
function at the critical point (that gives more information than the
critical exponent). In the case of a free product of finite groups
(say, ${\mathrm{PSL}}(2, \mathbb{Z})$), $p_c$ is obtained as a root
of an explicitly written polynomial (in particular, $p_c$ is an
algebraic number).

\subsection{Free products: critical probability}\label{here}

The Cayley graph of a free product of groups has a tree-graded
structure \cite{sapir:topology}: it is a union of subgraphs $M_j,
j\in {\mathbb N}$, each $M_j$ is a copy of the Cayley graph of one
of the group $G_i$, different $M_j$ and $M_k$ intersect by at most
one point, and every simple loop in the Cayley graph is in one of
the $M_i$. Note that the results of this paper can be generalized to
arbitrary transitive locally finite tree-graded graphs.

A non-trivial free product of groups (except $C_2*C_2$) has
infinitely many ends. This implies $p_u=1$ (as noticed for example
by Lyons \cite{lyons:book}). But the critical probabilities $p_c$
and $p_{\mathrm{exp}}$  were not known even in the case of the
modular group ${\mathrm{PSL}}(2, \mathbb{Z})$ which is the free
product of cyclic groups of orders $2$ and $3$.

We start with the following result giving an equation for $p_c$ in
the case of free product of two groups.

\begin{theorem}\label{fortwo} Let $G_1=\langle S_1\rangle$,
$G_2=\langle S_2\rangle$ be two finitely generated groups. Consider
the Cayley graph of the free product $G_1 * G_2$ with respect to the
generating set $S_1\cup S_2$.

Then the critical probability $0<p_c\le 1$ is the unique solution of
the following equation
\begin{align}\label{eq:two}
(\chi_1(p)-1)(\chi_2(p)-1)=1,
\end{align}
where $\chi_i(p)=\mathrm{E}_p(|C|_{G_i} )$ denotes the expected size of the
cluster containing the origin in the Cayley graph of group $G_i$ with
respect to the generating set $S_i$.
\end{theorem}

The critical probability $p_c$ of the free product is 1 if and only
if $|G_1|=|G_2|=2$ (in that case the group $G_1*G_2$ is virtually
cyclic).

The next theorem gives the expected size of a cluster for
subcritical $p$ in the Cayley graph of a free product of two groups.

\begin{theorem} \label{expected}
Let $G_1*G_2$ be as in Theorem \ref{fortwo}. Then for $p<p_c$, the mean cluster size satisfies
\begin{equation}\label{eq65}
\mathrm{E}_p(|C|_{G_1*G_2})=
\frac{\chi_1(p)\chi_2(p)}{\chi_1(p)+\chi_2(p)-\chi_1(p)\chi_2(p)}.
\end{equation}
\end{theorem}

Note that by Theorem \ref{fortwo}, the denominator in formula
(\ref{eq65}) is equal to 1 if $p=0$ and is decreasing to $0$ as
$p\to p_c(G_1*G_2)$.

The next two corollaries generalize the above theorems to the free
product of an arbitrary number of groups. They follow by
induction.

\begin{corollary}\label{main1}
Let $\mathcal{G}$ be a Cayley graph of the free product of $n$
non-trivial finitely generated groups $G_i$, $i=1,\dots,n$, with
respect to the set of generators $\bigcup_{i=1}^{n}S_i$, where $S_i$
is a generating set of $G_i$.

The expected size of the cluster at the origin is  equal to
\begin{align}\label{eq:ec}
\mathrm{E}_p(|C|_{G_1*\cdots*G_n})=\frac{\prod_{i=1}^n\chi_i(p)}{\sum_{j=1}^n\prod_{i=1,i\neq j}^n\chi_i(p)-(n-1)\prod_{i=1}^n\chi_i(p)}
\end{align}
for $p<p_c$ and it is infinity for $p\geq{p_c}$.
\end{corollary}

\begin{corollary} \label{main}
Let $\mathcal{G}$ be a Cayley graph of the free product of $n$
non-trivial finitely generated groups $G_i$, $i=1,\dots,n$, with
respect to the set of generators $\bigcup_{i=1}^{n}S_i$, where $S_i$
is a generating set of $G_i$. Then the critical probability
$0<p_c\leq1$ of $\mathcal{G}$ is the unique solution of the following
equation
\begin{equation}\label{eq:pc}
\sum_{j=1}^n\prod_{i=1,i\neq j}^n\chi_i(p)-(n-1)\prod_{i=1}^n\chi_i(p)=0,
\end{equation}
where $\chi_i(p)=\mathrm{E}_p(|C|_{G_i})$ denotes the expected size of the component containing origin in the Cayley graph of $G_i$ with respect to the set of generators $S_i$.
\end{corollary}

We give two proofs of Theorem \ref{fortwo}. The first proof is direct, and the second one uses the theory of branching processes \cite{parzen:book}. Theorem \ref{expected}  and Corollaries \ref{main1}, \ref{main} also can be obtained as applications of the theory of branching processes.

The next theorem shows, in particular, that in the case of free
products of residually finite groups, the critical probability can
be obtained as a limit of a fast converging sequence of algebraic
numbers (critical probabilities of free products of finite groups).

Suppose that a finitely generated group $G=\la S\ra$ has surjective
homomorphisms $\phi_i\colon G\to F_i$ such that $\phi_i$ is
injective on a ball of radius $i$ in the Cayley graph of $G$ (i.e.
on the set of all products of elements from $S\cup S^{-1}$ of length at most
$i$). Then we shall say that $G$ is {\em well approximated} by
$F_i$. By the Cayley graph of $F_i$ we shall always mean the Cayley
graph with respect to $\phi_i(S)$.

\begin{theorem}\label{limit} Suppose that each of the (non-trivial) finitely generated group
$G_i=\la S_i\ra$ is well approximated by groups $H_{i}^{j}$,
$i=1,2$. Then
$$
p_c(G_1*G_2)=\lim_{j\rightarrow \infty}p_c(H_1^j*H_2^j).
$$
More precisely there exist $C,\gamma>0$ such that
$$
0\leq p_c(H_1^j*H_2^j)-p_c(G_1*G_2)\leq C e^{-\gamma j}.
$$
\end{theorem}

A similar result holds (by induction) for a free product of any
finite number of groups.

Finally note that although inequalities (\ref{eq:est1}) and
(\ref{eq:est2}) become equalities for the free group (with free
generators) \cite{benjamini:beyond}, already for the free products
of finite cyclic groups both inequalities become strict (see
Proposition \ref{cheeger}). Thus these inequalities only give rough
estimates for $p_c$.

The next proposition gives the Schonmann's critical point for free
products.

\begin{proposition} \label{pexp}
Consider the free product $G_1*\cdots*G_n$ of $n$ finitely generated groups. Then \begin{align*}
p_{\mathrm{exp}}(G_1*\cdots*G_n)=\min_{1\leq i\leq n}\{p_{\mathrm{exp}}(G_i)\},
\end{align*}
where for finite $G_i$ we define $p_{\mathrm{exp}}(G_i)=1$.
\end{proposition}
In particular for the free product of finite groups $p_{\mathrm{exp}}$ is equal to one.

\begin{corollary}
Consider the free product $G_1*\cdots*G_n$ of $n$ nontrivial finitely generated groups, which is not virtually $\mathbb{Z}$. Then
\begin{align*}
p_c(G_1*\cdots*G_n)< \min_{1\leq i\leq n} \{p_c(G_i)\}\leq\min_{1\leq i\leq n} \{p_{\mathrm{exp}}(G_i)\}=p_{\mathrm{exp}}(G_1*\cdots*G_n)\leq p_u(G_1*\cdots*G_n)=1
\end{align*}
\end{corollary}
The first inequality will be shown to be strict for any free product, which is not virtually $\mathbb{Z}$. If it is virtually cyclic, i.e. $C_2*C_2$, all mentioned critical values are equal to one.

\subsection{Free products: critical exponents}

H\"aggstr\"om and Peres \cite{haggstrom:mono} proved that the
function $\theta(p)$ is continuous for $p>p_c$ for any Cayley graph. A free product of nontrivial groups (which is not virtually cyclic) is non-amenable, so the percolation dies at $p_c$ (for the proof see \cite{benjamini:nonamen}) and therefore the percolation function $\theta(p)$ is continuous  also at $p_c$. This allows us to consider the critical exponents.

\comment of free products.
We obtain
\begin{proposition} \label{expon} Consider the free product $G_1*\cdots*G_n$ of $n$
 (nontrivial) finitely generated groups. Assume $n>2$ or $|G_1|>2$ or $|G_2|>2$.
Then the critical exponent $\gamma$ is equal to $-1$. If moreover
$G_1,\dots,G_n$ are finite, then the critical exponent $\beta$ is
equal to $1$. The derivative of the percolation function at $p_c$
from the right is a rational function of the first and second moment
of the cluster sizes in $G_i$'s.
\end{proposition}
\endcomment

We introduced two critical exponents $\beta$ and $\gamma$, but there are several others describing the behavior near the critical point $p_c$. Physicists believe that the numerical values of critical exponents depend only on the underlying space and not on the structure of the particular lattice.

Values of critical exponents are well known for trees ($\beta=1$, $\gamma=-1$), so called {\em mean-field values}. Hara and Slade \cite{hara:meanfield} proved that the critical exponents take their mean-field values in $\mathbb{Z}^d$ for $d>19$ by verifying the triangle condition introduced by Aizenman and Newman \cite{newman:triangle}. Schonmann \cite{schonmann:planar} proved that the critical exponents take their mean-field values for all non-amenable planar graphs with one end, and for unimodular graphs with infinitely many ends. The later case covers Cayley graphs of free products.

The equation for percolation function found in the proof of Theorem
\ref{fortwo} allows us to evaluate the right derivative of the percolation
function using implicit differentiation. One can also compute the  left
derivative of the expected cluster size function of a free product at $p_c$.
The formulas for these derivatives immediately imply Schonmann's result that
the values of critical exponents of a free product
are mean-field. In particular the formula for the derivative of the cluster
size of a free product is the following. It involves the derivative of
cluster sizes in the factor groups:

\begin{align*}
\left.\frac{d}{d _-p}(\mathrm{E}_{p_c}(|C|_{G_1*G_2}))^{-1}\right|_{p=p_c}
=\frac{\left.\frac{d}{d _-p} \chi_1(p) \right|_{p=p_c}(1-\chi_2(p_c))+\left.\frac{d}{d _-p} \chi_2(p) \right|_{p=p_c}(1-\chi_1(p_c))   }{\chi_1(p_c)\chi_2(p_c)}
\end{align*}

\comment
 We included Proposition \ref{expon} because the proof is direct and involves exact computation of the limiting behavior.
 In order to find the critical exponent we evaluate the first term in the Taylor expansion of $\theta(p)$ at $p_c$, and similarly the whole Taylor expansion can be found.
\endcomment

Schonmann \cite{schonmann:nonamen} also proved the mean-field criticality for highly non-amenable graphs, i.e. such that $h(G)> D(G)/\sqrt{2}$, where $D(G)$ is the maximal degree of vertices in $G$ (in case of a Cayley graph it is just the degree of the origin $\deg(o)$). The question for general non-amenable graph remains open. Define the {\it spectral radius} of $G$ by
$$R(G)=\limsup_{n\to \infty}(\#\text{ of closed walks of length }n\text{ at }o)^{1/n}.
$$
Schonmann used in his argument that if $p_c<1/R(G)$ then the
triangle condition is satisfied and thus the critical exponents take their
mean-field values. For any non-amenable group, Pak and Smirnova-Nagnibeda
 construct a finite generating set   with
the property $p_c R(G)<1$. Combining these two results we get the
following statement.

\begin{proposition}
Every finitely generated non-amenable group has a finite generating set such that the Cayley graph with respect to this generating set has mean-field valued critical exponents.
\end{proposition}

 Sapir conjectured that the mean-field criticality should be true
for a class of hyperbolic groups which includes free products (see Conjecture
\ref{conj} below). Recall that a group $G$ is called
{\em Gromov-hyperbolic} if for some $\delta>0$, in every geodesic triangle
of the Cayley graph, one side is in the $\delta$-neighborhood of the
union of two other sides. A group is called {\em elementary} if it
contains a cyclic subgroup of finite index. For example, free
products of finite groups are hyperbolic and if the groups
are of order $>2$, the free products are not elementary.

A well known conjecture in percolation theory claims that the
critical exponents of all  lattices in ${\mathbb R}^2$ are the same
($\beta=5/36$, $\gamma=-43/18$) as proved for triangular lattices by Smirnov and
Werner \cite{werner:critexp}). The motivation for this conjecture is
that for every lattice $L$ in ${\mathbb R}^2$ the Gromov-Hausdorff
limit of rescaled copies of $L$, $L/2$, $L/3$, \dots is isometric to
${\mathbb R}^2$ with the $L_1$-metric, i.e. any two lattices in
${\mathbb R}^2$ have the same asymptotic cones
(\cite{gromov:asymptotic}). Sapir conjectured that this is true
in general: if two groups have isometric asymptotic cones then their
critical exponents should coincide. In particular, since all
asymptotic cones of all non-elementary Gromov-hyperbolic groups are
isometric (they are isometric to the universal ${\mathbb R}$-tree of
degree continuum by a result of Dyubina and Polterovich
\cite{dyubina}), the following statement should follow:

\begin{conjecture}[M. Sapir]\label{c1} If the asymptotic cones of Cayley graphs of two groups $G$, $G'$ are isometric, then their critical exponents are the same.

In particular every Cayley graph of a non-elementary hyperbolic group has mean-field valued critical exponents.
\label{conj}
\end{conjecture}

Note that Benjamini and Schramm \cite{benjamini:beyond} asked the question whether all Cayley graphs of non-amenable groups have mean-field criticality. The positive answer to this question would of course imply the second part of Conjecture \ref{c1}.

The author would like to thank Mark Sapir for valuable discussions and careful correcting of the text. I am also grateful to Tatiana Smirnova-Nagnibeda for useful comments and Russel Lyons for suggesting the approach of branching processes.

\section{Critical probability $p_c$}
\subsection{A recursive expression for percolation function}\label{recur}

First let us describe the structure of a Cayley graph $\mathcal{G}$ of the free
product of two groups. Every vertex $v$ is contained in exactly two
basic subgraphs induced by vertices: the basic subgraph of the first type,
$\mathcal{G}_1$, is an induced subgraph of $\mathcal{G}$ consisting
of vertices obtained from $v$ by multiplying on the right by
elements of $G_1$, and the basic subgraph of type two,
$\mathcal{G}_2$, is defined similarly for $G_2$. Each of the basic
subgraphs of type one (type two) is isomorphic to the Cayley graph
of $G_1$ (resp. $G_2$). Every vertex of the Cayley graph of $G$ is
the (unique) common vertex of a subgraph of type $1$ and a subgraph
of type $2$. In the case of finite cyclic groups with single
generator, $\mathcal{G}_i$'s are cycles.

Since the Cayley graph of a finitely generated group is locally
finite, the origin is in an infinite cluster if and only if there
exists an infinite open simple path starting from the origin. Note that
if $G_1$ and $G_2$ are finite, any infinite simple path in the
Cayley graph of $G$ has to intersect infinitely many basic subgraphs
and if it leaves a subgraph, it can never return to it. In the general case, any simple infinite path starting at the vertex $v$ is of one of the following two types: type one has first edge in the basic subgraph $\mathcal{G}_1$ (and uses no edges of the basic subgraph $\mathcal{G}_2$ which contains $v$); type two starts with an edge from basic subgraph $\mathcal{G}_2$.

Let $A$ be the probability that there is an infinite path of type one starting at the origin $o$
%and visiting an edge of $\mathcal{G}_1$ at $o$. \
and let $B$ be defined similarly for the type two paths.
These two events are independent. Moreover we can recursively express $A$
using $B$ and vice versa as follows.
\begin{align}\label{eq:all}
1-A=\sum_{\omega \text{ subgraph of } \mathcal{G}_1}
\mathrm{P}_p(\omega) (1-B)^{|C|_{G_1}-1},
\end{align}
where $\mathrm{P}_p(\omega)$ is the probability of $\omega=(V(\omega),E(\omega))$ to be the open subgraph of
$\mathcal{G}_1=(V(\mathcal{G}_1),E(\mathcal{G}_1))$.
The left side is the probability that there is no infinite path (starting at
the origin) of type one.
The right side is obtained by conditioning on at which vertex the infinite path could leave $\mathcal{G}_1$.
Specifically for a fixed realization $\omega$  we compute probability that   there is no infinite path of type two starting at a vertex of the connected component of the origin in $\mathcal{G}_1$.
Probability of the existence of an infinite path of type two starting at any specific vertex is equal to $B$ because Cayley graphs are homogeneous (and so it does not matter which vertex we take as the origin).
Therefore the right side is equal to the probability that there is no simple infinite path starting at the origin and leaving $\mathcal{G}_1$ at any vertex (different from the origin) connected to the origin.

Note that for finite groups
$$
\mathrm{P}_p(\omega)=p^{|E(\omega)|}(1-p)^{(|E(\mathcal{G}_1)|-|E(\omega)|)},
$$
and $|C|_{G_1}$ is the number of vertices in the connected component of
$\omega$ containing $o$ in $\mathcal{G}_1$.

For infinite $G_1$ the summation should be replaced by integration, or we shall use  short expression using the expectation in the probability space restricted to $\mathcal{G}_1$:
\begin{align*}
1-A=\mathrm{E}_p[(1-B)^{|C|_{G_1}-1}].
\end{align*}
One can modify the equation (\ref{eq:all}) by summing only over all
connected subgraphs of $\mathcal{G}_1$ containing $o$:
\begin{align*}
1-A=\sum_{K \text{ finite connected subgraph of } \mathcal{G}_1
\text{ containing } o} \mathrm{P}_p'(K) (1-B)^{|C|_{G_1}-1}.
\end{align*}
Here $\mathrm{P}_p'(K)=p^{|E(K)|}(1-B)^{|\partial
K|}$ where $|\partial K|$ denotes the size of the
(external) boundary of $K$ in $\mathcal{G}_1$, i.e. the number
of edges of $\mathcal{G}_1$ not in $K$ with at least one end
vertex (maybe both) in $K$.

Similar equalities hold for $\mathcal{G}_2$. Let us now formally rewrite
the summation to be over the size of the connected component
containing $o$ of a random subgraph $\omega$ of $\mathcal{G}_i$. Denote by $Q_i(n)$ the probability that this component is of size $n$, i.e. $Q_i(n)=\mathrm{P}_p(|C|_{G_i}=n)$.
Define a recurrent {\it walk through} function $g_i$, $i=1,2$ for $0\leq p,t\leq1$ by:
\begin{align}\label{eq:gi}
g_i(p,t)=1-\sum_{n=1}^{|G_i|}(1-t)^{n-1}Q_i(n)
\end{align}
Notice that $g_i$ is very close to the moment generating function of the cluster size in $\mathcal{G}_i$ (see for example \cite{feller:book}). This will become handy when we take the derivatives at $t=0$ in Proposition \ref{properties}.

If $G_i$ is an infinite group then we sum over infinitely many
(non-negative) values, but the sum is always bounded by 1 from above
(since $\sum_{n=1}^{\infty}Q_i(n)\leq 1$). On
the other hand $g_i$ is always bigger than the probability of having
an infinite cluster at the origin just in $G_i$ (i.e.
$\theta_{G_{i}}(p)$).

In the case when $G_i$ is a cyclic group of finite order
$m$, and so $\mathcal{G}_i$ is a finite cycle, connected
subgraphs containing $o$ are just arcs. Probability of a specific arc of length $n<m$ is $(1-p)^2p^{n-1}$ and the number of those containing origin is just $n$. Therefore
\begin{align*}
Q_i(n)&=n(1-p)^2p^{n-1} \text{ for }n<m,\\
Q_i(m)&=(m(1-p)+p)p^{m-1}.
\end{align*}
It simplifies the summation as follows.
\begin{align}
\label{eq:cyc}
g_i(p,t)=1-\sum_{j=1}^{m-1} (j(1-p)^2(p(1-t))^{j-1})  -
(m(1-p)+p)(p(1-t))^{m-1}
\end{align}

The percolation function is given by
\begin{align}\label{eq:theta}
\theta&=A+B-AB,\notag\\ \intertext{where}
A&=g_1(p,B),\notag\\
B&=g_2(p,g_1(p,B)).
\end{align}
It remains to determine for which
values of $p$ the last equation in (\ref{eq:theta}) has a positive solution.

\subsection{The critical probability}

We start with describing general properties of the walk-through function $g_i$.
In particular we will show that the function $\varrho_p(B)=g_2(p,g_1(p,B))-B$ is concave and equal to $0$ for $B=0$. This allows us to decide whether the equation (\ref{eq:theta}) has a positive solution based on the derivative of $\varrho_p(B)$ at zero.

We will use the following result of Aizenman and Barsky
\cite{aizenman}. In fact they proved it for $\mathbb{Z}^d$ but their
argument works for any transitive graph as noticed by Lyons and
Peres \cite{lyons:book}.

\begin{lemma}[Aizenman, Barsky]\label{lemma}
If $p<p_c$ then the mean cluster size is finite in any transitive graph.
\end{lemma}

\begin{proposition}[Properties of $g_i$]\label{properties}
Let $g_i$, $i=1,2$ be defined as above for arbitrary finitely
generated groups $G_i$ (with some choice of generators). Assuming
$0< t < 1$, $0 < p < 1$, we have:
\begin{enumerate}
\item $0\leq g_i(p,t)\leq 1$.

\item $\left.\frac{\partial^k g_i(p,t)}{\partial t^k}\right|_{t=0}=(-1)^{k+1} \mathrm{E}_p[(|C|_{G_i}-1)(|C|_{G_i}-2)\dots(|C|_{G_i}-k)]$ for $k=1,2,..$ and $p<p_c(G_i)$.

\item $g_2(p,g_1(p,t))-t$ is concave in $t$ for $p<min(p_c(G_1),p_c(G_2))$.
\end{enumerate}
\end{proposition}
\begin{proof}
Using the definition of $g_i$ given by (\ref{eq:gi}) we may evaluate it for $t=0$ and $t=1$. Recall that $Q_i(n)$ is the probability that the component containing origin in $\mathcal{G}_i$ has size $n$,
\begin{align}
g_i(p,0)&=
1-\sum_{n=1}^{|G_i|}Q_i(n)=Q_i(\infty)\geq 0
\notag
\\
g_i(p,1)&=
1-Q_i(1)=1-(1-p)^{D(G_i)}\leq 1
\notag
\end{align}
where $D(G_i)$ is the degree of the origin (and any other vertex) in the graph  $\mathcal{G}_i$. Now let us take the formal derivative.
\begin{align}
\frac{\partial g_i(p,t)}{\partial t} &=
\sum_{n=1}^{|G_i|}(n-1)(1-t)^{n-2}Q_i(n)\geq 0
\notag\\
\frac{\partial g_i(p,t)}{\partial t} &\leq \left.\frac{\partial g_i(p,t)}{\partial t}\right|_{t=0} = \sum_{n=1}^{|G_i|}(n-1) Q_i(n)=
\chi_i(p)-1),
\end{align}
Note that $\chi_i(p)=|G_i|$ and for $p<p_c(G_i)$, $\chi_i(p)<\infty$ by Lemma \ref{lemma}. This in particular implies the continuity and differentiability of $g_i$
and $\chi_i(p)$ for $p<p_c$.

From above we can conclude that $0\leq g_i\leq 1$ as expected from the fact that for
particular values of $t$, $g_i$ gives the probability of an infinite
path in a part of the Cayley graph as defined above.

Now we compute the $n$-th derivative (term by term) and evaluate it
at $t=0$.
\begin{align}
\frac{\partial^k g_i(p,t)}{\partial t^k} &= (-1)^{k+1}
\sum_{n=k+1}^{|G_i|}(n-1)(n-2)\dots(n-k)(1-t)^{n-k-1}Q_i(n)
\label{eq:der}
\end{align}
If $p<p_c(G_i)$ then $Q_i(\infty)=0$ and we have (by the definition of the expectation)
\begin{align*}
\left.\frac{\partial^k g_i(p,t)}{\partial t^k}\right|_{t=0}&=(-1)^{k+1} \mathrm{E}_p[(|C|_{G_i}-1)(|C|_{G_i}-2)\dots(|C|_{G_i}-k)]
\end{align*}
This proves part ii.

Using the above formula for the derivative of $g_i$ we see that the
odd numbered derivatives are positive and the even numbered
derivatives are negative. Thus in particular $g_i$ is increasing and
concave in $t$ for all $p$.

%Formula (\ref{eq:der}) gives the derivative of $g_i$ as %$\mathrm{E}_p(|C|_{G_i}-1)$ if $\mathrm{P}_p(|C|_{G_i}=\infty)=0$, that is for %$p<p_c(G_i)$ if $G_i$ is infinite and for all $p$ if $G_i$ is finite (we set
%$p_c(G_i)=1$ for finite $G_i$).

As soon as $|G_i|>2$, $g_i$ is strictly concave and
$\chi_i(p)$ is increasing. Since the composition of
two concave functions is concave, $g_2(p,g_1(p,t))-t$ is concave in
$t$ for $p< \min(p_c(G_1),p_c(G_2))$, and if
$(|G_1|-1)(|G_2|-1)>1$ then it is strictly concave .
\end{proof}

Now we are ready to show for which $p$ the percolation function is
positive, that is to find $p_c$ of free product of two
groups.

\begin{proof}[Proof of Theorem \ref{fortwo}]
Clearly $p_c(G_1*G_2)\leq \min(p_c(G_1),p_c(G_2))$ so we need to
decide for which $p<\min(p_c(G_1),p_c(G_2))$ the equation
$g_2(p,g_1(p,t))-t=0$ has a positive solution.

In what follows we always assume $p<\min(p_c(G_1),p_c(G_2))$ and set $\varrho_p(t)=g_2(p,g_1(p,t))-t$.

Since $g_i(p,0)=Q_i(\infty)=0$ we have $\varrho_p(0)=g_2(p,g_1(p,0))=0$. On the other hand $g_2(p,g_1(p,1))\leq 1$
and thus $\varrho_p(1)\leq 0$. By Proposition \ref{properties} iii. $\varrho_p(t)$ is concave in $t$. Therefore there is at most one change in the monotonicity of $\varrho_p(t)$ on the unit interval, in particular the function is either decreasing all the time or there is $t_0$  such that  $\varrho_p(t)$ is increasing the interval $(0,t_0)$ and decreasing on $(t_0,1)$, see Picture \ref{fig:fce}.  Thus
the equation  $\varrho_p(t)=0$ has a positive solution ($0$ is always a solution) if and only if the function $\varrho_p(t)$ has positive derivative at $t=0$.

\begin{figure}[t]
\begin{center}
  \includegraphics[trim=0 180 180 0,scale=.4,angle=270]{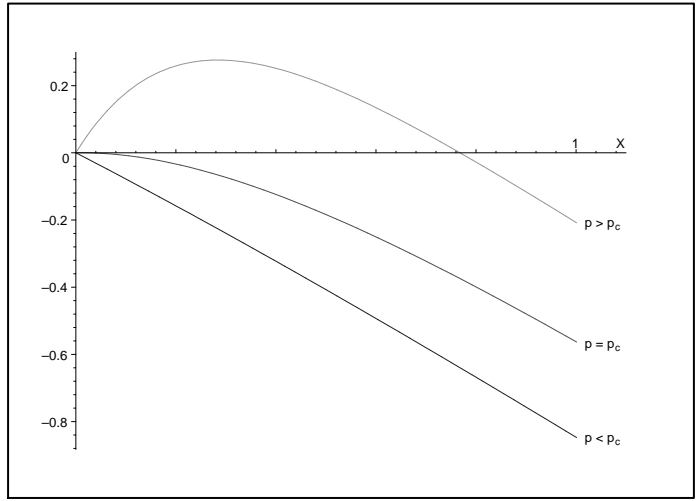}
  \caption{Graph of a  $g_2(p,g_1(p,x))-x$ in case $G=C_3*C_5$ for different values of $p$.}
  \label{fig:fce}
\end{center}
\end{figure}

Now using Proposition \ref{properties} ii. we have:
\begin{align*}
\left.\frac{\varrho_p(t)}{\partial t}\right|_{t=0}&=
\left.\frac{\partial g_2(p,g_1(p,t))-t}{\partial t}\right|_{t=0}
\\&=
\left.\frac{\partial g_2(p,t)}{\partial t}\right|_{t=g_1(p,0)=0}
\left.\frac{\partial g_1(p,t)}{\partial t}\right|_{t=0}-1
\\&= (\chi_2(p)-1)(\chi_1(p)-1)-1
\end{align*}
where $\chi_i(p)$ is the expected size of the
component of $\mathcal{G}_i$ containing the origin. It is an
increasing function of $p$.

Therefore the value of derivative of $\varrho_p(t)$ at $t=0$ is increasing function of $p$. For $p=0$ it equals $-1$ and for $p_1=\min(p_c(G_1),p_c(G_2))$ we have
\begin{align*}
\left.\frac{\varrho_{p_1}(t)}{\partial t}\right|_{t=0}=
(|G_1|-1)(|G_2|-1)-1\geq 0
\end{align*}
Therefore there exists unique $0<p_0\leq \min(p_c(G_1),p_c(G_2))$
such that
\begin{align*}
\left.\frac{\varrho_{p_0}(t)}{\partial t}\right|_{t=0}=0.
\end{align*}
or equivalently
\begin{align*}
(\chi_1(p_0)-1)(\chi_2(p_0)-1)=1.
\end{align*}

The case $C_2*C_2$ (the group is virtually $\mathbb{Z}$) is the only
one when $p_0=1$ and its percolation function $\theta(p)=0$ for
$0\leq p<1$, and $\theta(1)=1$. In all other cases we conclude that
the equation $\varrho_{p}(t)=0$ has unique nonzero solutions if
and only if $p > p_0$. So $p_0$ is the critical probability $p_c(G_1*G_2)$ of the Cayley graph of free product $G_1*G_2$. Note that if $(|G_1|-1)(|G_2|-1)\geq 1$ (i.e. both groups are nontrivial and one of them has more than $2$ elements) then $p_c(G_1*G_2)<\min{p_c\{G_i\}}$.

This finishes the proof of Theorem \ref{fortwo}.
\end{proof}

%We shall need the following statement when we consider critical
%exponents.
%\begin{proposition}\label{conti} The function $\theta(p)$ is
%continuous at $p_c$.
%\end{proposition}
%\proof
% Free product of nontrivial groups is non-amenable, so the
%percolation dies at $p_c$, i.e. $\theta(p_c)$=0 (for the proof see
%\cite{benjamini:nonamen}). We will show that the limit of $\theta$
%at $p_c$ from the right is equal to $0$.
%
% By (\ref{eq:theta}), the percolation
%function $\theta$ is a continuous function of $B$, (the probability
%that the infinite path in the cluster around $o$ exists starting
%with an edge from a given part of the Cayley graph). As a function
%of $p$, $B$ is nondecreasing and bounded. Thus the limit of $B$ from
%the right at $p_c$ exists. Denote $B_0$ the limit of $B=t(p)$ as $p$
%approaches $p_c$ from the right and assume $B_0>0$ then
%$$
%g_2(p_c,g_1(p_c,B_0))-B_0=\lim_{p\rightarrow p_{c+}}
%g_2(p,g_1(p,t(p)))-t(p)=0
%$$
%Thus $B_0$ would be a positive solution of $g_2(p,g_1(p,t))-t=0$ at
%$p_c$ contradicting the fact that percolation dies at $p_c$.
%Therefore $\lim_{p\rightarrow
%p_{c+}}\theta(p)=\theta(p_c)=0$.\endproof
%
%Similarly $\lim_{p\rightarrow 1_-}\theta(p)=\theta(1)=1$.

\subsection{The expected cluster size}\label{ctyri}

The expected cluster size function can be expressed recurrently in a
similar way as the percolation function $\theta(p)$. In this
section, we consider the free product of arbitrary number of groups.

\begin{proposition} Let $G_i=\la S_i\ra$ be a finitely generated group,
$i=1\dots n$. Denote by $\chi_i(p)$ the expected size of a
cluster at the origin in the Cayley graph of $G_i$ with respect to
the generating set $S_i$. Then the expected size of the cluster at
the origin in a Cayley graph of the free product $G_1*\cdots*G_n$ with
respect to the generating set $S_1\cup \dots \cup S_n$ satisfies
\begin{align}
\mathrm{E}_p(|C|_{G_1*\cdots*G_n})=\frac{\prod_{i=1}^n\chi_i(p)}{\sum_{j=1}^n\prod_{i=1,i\neq j}^n\chi_i(p)-(n-1)\prod_{i=1}^n\chi_i(p)} \label{eq:exp}
\end{align}
for $p<p_c(G_1*\cdots*G_n)$.
\label{expsize}
\end{proposition}

\begin{proof}
First assume n=2. If we remove the origin from the Cayley graph it
splits into two parts: part $P_1$ consists of vertices $v$ such that
any simple path from the origin to $v$ first visits a vertex
corresponding to an element of $G_1$ (in fact $S_1$), part $P_2$ is
defined similarly for $G_2$. Denote $|C|_{P_i}=|C\cap
P_i|_{G_1*G_2}$, size of the cluster containing the origin intersected with the part $P_i$. Now we can write
\begin{align*}
\mathrm{E}_p(|C|_{G_1*G_2})&=\, \mathrm{E}_p(|C|_{P_1})+\mathrm{E}_p(|C|_{P_2})+1,
\end{align*}
where the constant $1$ represents the origin. One can represent
$|C|_{P_1}$ as a random sum of random variables:
$$
|C|_{P_1}=\sum_{o\neq x\in C\cap G_1}Y_x,
$$
where each $Y_x$ has the same
distribution as $|C|_{P_2}+1$. Using Wald's identity (see for example \cite{feller:book}) we have:
\begin{align*}
\mathrm{E}_p(|C|_{P_1})&=\, (\chi_1(p)-1)(\mathrm{E}_p(|C|_{P_2})+1)\\
&=\, \chi_1(p)-1+(\chi_1(p)-1)(\chi_2(p)-1)(\mathrm{E}_p(|C|_{P_1})+1)
\end{align*}
A similar equality holds for $\mathrm{E}_p(|C|_{P_2})$. Combining
into one equation and solving for $\mathrm{E}_p(|C|_{G_1*G_2})$ we
obtain:
\begin{align}
\mathrm{E}_p(|C|_{G_1*G_2})&=\, \chi_1(p)+\chi_2(p)-1+(\chi_1(p)-1)(\chi_2(p)-1)(\mathrm{E}_p(|C|_{G_1*G_2})+1)\notag
\\
&=
\frac{\chi_1(p)\chi_2(p)}{\chi_1(p)+\chi_2(p)-\chi_1(p)\chi_2(p)} \label{eq:exp2}
\end{align}
Note that the formula in the denominator coincides with the one in
the equation (\ref{eq:two})  after rearrangement. Therefore it is positive for $p<p_c$
and equal to $0$ at $p_c$. Thus the expected cluster size is given by formula (\ref{eq:exp2}) for $p<p_c$ and it tends to infinity as $p$ approaches
$p_c$ (the expected size of the cluster is equal to infinity for
$p\geq p_c$).

Suppose now that the result holds for any free product of at most
$n$ groups. Consider $G_1*\cdots*G_{n+1}$. We can view it as a free
product of two groups  $G_1*\cdots*G_{n}$ and $G_{n+1}$. By induction
we know the expected size of the component in both of those groups
and we can apply formula (\ref{eq:exp2}). Therefore
\begin{align*}
\mathrm{E}_p(|C|_{G_1*\cdots*G_{n+1}})&=\frac{\mathrm{E}_p(|C|_{G_1*\cdots*G_{n}})\chi_{n+1}(p)}
{\mathrm{E}_p(|C|_{G_1*\cdots*G_{n}})+\chi_{n+1}(p)-\mathrm{E}_p(|C|_{G_1*\cdots*G_{n}})\chi_{n+1}(p)}\\
&=\prod_{i=1}^n\chi_i(p) \, \chi_{n+1}(p)\delta^{-1},\\
\end{align*}
%\intertext{
where
%}\qquad
\begin{align*}
\delta&=
\prod_{i=1}^n\chi_i(p)+
\chi_{n+1}(p)\left(\sum_{j=1}^n\prod_{i=1, i\neq j}^n\chi_i(p)-(n-1)\prod_{i=1}^n\chi_i(p)\right)
%\\&\quad
-\chi_{n+1}(p)\prod_{i=1}^n\chi_i(p)
\\
&=
\sum_{j=1}^{n+1}\prod_{i=1, i\neq j}^{n+1}\chi_i(p)
-n\prod_{i=1}^{n+1}\chi_i(p).
\end{align*}
For $p<p_c$ all expressions are finite and the denominator is non-zero (follows from the discussion of the free product of two groups).

\end{proof}

Corollary \ref{main} about $p_c$ for $G_1*\cdots*G_{n}$ now follows
from Proposition \ref{expsize}. Let us look at the formula for the
expected cluster size more closely. The numerator in formula
(\ref{eq:exp}) is finite, positive and increasing for
$p<\min\{p_c(G_i)\}$. The denominator is equal to $1$ for $p=0$ and
is decreasing in $p$ for $p<p_c$ since each
$\chi_i(p)$ is increasing in $p$ and
\begin{align*}
\frac{\partial \left(\sum_{j=1}^n\prod_{i=1,i\neq j}^n\chi_i(p)-(n-1)\prod_{i=1}^n\chi_i(p)\right)
 }{\partial \chi_k(p)}&=
 \sum_{j=1,j\neq k}^n (1-\chi_j(p) )\prod_{i=1,i\neq j,k}^n\chi_i(p)
 \\&<0.
\end{align*}
As $p$ approaches $\min\{p_c(G_i)\}$ the expression in the denominator becomes negative, and thus there exists unique $p$ such that
$$\sum_{j=1}^n\prod_{i=1,i\neq j}^n\chi_i(p)-(n-1)\prod_{i=1}^n\chi_i(p)=0.
$$
The solution is the critical probability $p_c$ which proves
Corollary \ref{main}.

\subsection{Comparison with branching processes}

The expected cluster size can be viewed as a population size of a branching process in the following way.
Consider the free product of two groups $G_1*G_2$. Let us define a branching process. The origin is the only element of generation zero.
The first generation consists of those vertices of the basic subgraph $G_1\setminus\{o\}$ that are connected to the origin by open paths,
the $2i$-th (resp. $2i+1$-st) generation contains vertices,
which are connected to some vertex of the previous generation by a nontrivial open path included in a subgraph of type $G_2$ (type $G_1$, resp.).
If we consider only even generations, then we obtain a simple Galton-Watson branching process.  Again using the Wald's identity (see for example \cite{feller:book}) we obtain that the expected size of one generation is precisely $(\chi_1(p)-1)(\chi_2(p)-1)$.
By the Basic theorem of branching processes (see for example theorem 2A on page 201 in \cite{parzen:book}), such a branching process terminates if the expected generation size is at most one, otherwise the population size grows to the infinity.
Therefore the critical value of $p$ occurs for  $(\chi_1(p)-1)(\chi_2(p)-1)=1$ as claimed by Theorem \ref{fortwo}.
The result can be generalized to the free product of arbitrary number of groups or to a general tree-graded vertex-transitive graph.

\section{The critical point $p_{\mathrm{exp}}$}

Next we would like to find the $p_{\mathrm{exp}}$ of the free product. That is to decide for which $p$ the connectivity function decays exponentially with the distance.

\begin{proof}[Proof of Proposition \ref{pexp}]

Clearly $p_{\mathrm{exp}}(G_1*\cdots*G_n)\leq \min\{p_{\mathrm{exp}}(G_i)\}$. For the converse inequality suppose that $p<\min\{p_{\mathrm{exp}}(G_i)\}$. Then there exist $C_i, \gamma_i>0$, for $1\leq i\leq n$ such that
\begin{align*}
\mathrm{P}_p(x \leftrightarrow y\text{ in }G_i)\leq C_i e^{-\gamma_i\mathrm{dist}(x,y)}.
\end{align*}
Clearly if $G_i$ is finite this estimate can be done for any $p$.
Define $C=\max\{C_i\}$, $\gamma=\min\{\gamma_i\}$.
Any (simple) path from $x$ to $y$ can be divided into finite number of pieces, where each piece lies in one special subgraph $\mathcal{G}_i$ and two consequent subgraphs share exactly one vertex. Thus there is a finite set of {\it cut points} $z_1,\dots,z_k$, the same for all paths connecting $x$ and $y$.
Let $z_0=x$ and $z_{k+1}=y$, we have
\begin{align*}
\mathrm{P}_p(x\leftrightarrow y)&=\prod_{i=0}^{k}\mathrm{P}_p(z_i\leftrightarrow z_{i+1}),\\
\mathrm{dist}(x,y)&=\sum_{i=0}^{k}\mathrm{dist}(z_i,z_{i+1}).
\end{align*}
Now there exists $K>0$ and $0<\alpha\leq \gamma$ such that if $\mathrm{dist(z_i,z_{i+1})}>K$ then
\begin{align*}
\mathrm{P}_p(z_i\leftrightarrow z_{i+1})\leq Ce^{-\gamma\mathrm{dist}(z_i,z_{i+1})}\leq e^{-\alpha\mathrm{dist}(z_i,z_{i+1})}.
\end{align*}
Considering only the state of edges adjacent to $z_i$  we obtain following rough estimate:
\begin{align*}
\mathrm{P}_p(z_i\leftrightarrow z_{i+1})\leq 1-(1-p)^{\deg(o)}.
\end{align*}
The right side is strictly less than one thus there exists $0<\beta\leq \alpha$ such that
\begin{align*}
 1-(1-p)^{\deg(o)}\leq e^{-\beta K}.
\end{align*}
Therefore if $\mathrm{dist}(z_i,z_{i+1})\leq K$ then
\begin{align*}
\mathrm{P}_p(z_i\leftrightarrow z_{i+1})\leq e^{-\beta K}\leq e^{-\beta\mathrm{dist}(z_i,z_{i+1})},
\end{align*}
and for $\mathrm{dist}(z_i,z_{i+1})> K$ we have
\begin{align*}
\mathrm{P}_p(z_i\leftrightarrow z_{i+1})\leq e^{-\alpha\mathrm{dist}(z_i,z_{i+1})}\leq e^{-\beta\mathrm{dist}(z_i,z_{i+1})}.
\end{align*}
Combining the above two estimates we obtain
\begin{align*}
\mathrm{P}_p(x\leftrightarrow y)&=e^{-\beta\mathrm{dist}(x,y)}.
\end{align*}
Therefore the connectivity function has an exponential decay at $p$ and we have proved $$p_{\mathrm{exp}}(G_1*\cdots*G_n)\geq \min\{p_{\mathrm{exp}}(G_i)\}.$$
\end{proof}

\section{Approxiation results}

\begin{proof}[Proof of Theorem \ref{limit}]
Recall that we consider a sequence of factor groups
$H_i^{j}=G_i/N_i^j$ and want to show that
$p_c(G_1*G_2)=\lim_{j\rightarrow \infty}p_c(H_1^j*H_2^j)$. Let
$\phi$ be the factor map $G_i \rightarrow H_i^j$. Consider a Cayley
graph of $G_i$ with respect to generating set $S$ and a Cayley graph
of $H_i^j$ with respect to generating set $\phi(S)$. Then any path
in $H_i^j$ from origin $o$ to a vertex $x$ can be lifted to a unique
path in $G_i$ from the origin $o$ to a vertex $y$, s.t. $\phi(y)=x$.
On the other hand the image under $\phi$ of a simple path does not
have to be simple. (A similar argument was used by Campanino
\cite{campanino}.) Thus
\begin{align*}
\mathrm{E}_p(|C|_{H_i^j})&=\sum_{x\in H_i^j}
\mathrm{P}_p(o \leftrightarrow x)\\
&=\sum_{x\in H_i^j}\mathrm{P}_p(\text{at least
one path }o \leftrightarrow x\text{ is open})\\
&\leq \sum_{x\in H_i^j}\sum_{\{y: \phi(y)=x\}}
\mathrm{P}_p(\text{at least one path }o \leftrightarrow y\text{ is open})\\
&=\sum_{y\in G_i}\mathrm{P}_p(o \leftrightarrow y)\\
&=\chi_i(p)
\end{align*}
From Theorem \ref{fortwo} we know that $p_c$ is a solution of
equation (\ref{eq:two}). Therefore $p_c(G_1*G_2)\leq
p_c(H_1^j*H_2^j)$ for all $j$ (that also follows from
\cite{campanino}).

Now assume that $p$ is such that $\chi_i(p)<\infty$,
$i=1,2$. Then by Schonmann \cite{schonmann:nonamen} there exist $C, \gamma>0$ such that $P_p(|C|_{G_i}\geq n)\leq C e^{-\gamma n}$ and
\begin{align*}
0\leq\chi_i(p)-\mathrm{E}_p(|C|_{H_i^j})&\leq
2\sum_{k=j}^{\infty}k\mathrm{P}_p(|C|_{G_i}=k)
% \stackrel{\scriptscriptstyle j\to\infty}{\longrightarrow}0
\\
&\leq 2C(j e^{-\gamma j}+\sum_{k=j+1}^{\infty}e^{-\gamma k})
\leq \frac{4je^{-\gamma j}}{1-e^{-\gamma}}\leq C'e^{-\gamma'j},
\end{align*}
where $C'>0$ (later also $C''>0$). Then
\begin{multline*}
0\leq
(\chi_1(p)-1)(\chi_2(p)-1)-
(\mathrm{E}_p(|C|_{H_1^j})-1)(\mathrm{E}_p(|C|_{H_2^j})-1)
=\\
=(\chi_1(p)-1)(\chi_2(p)-\mathrm{E}_p(|C|_{H_2^j}))
+(\chi_1(p)-\mathrm{E}_p(|C|_{H_1^j}))(\mathrm{E}_p(|C|_{H_2^j})-1)
\leq\\
\leq(\chi_1(p)-1)(\chi_2(p)-\mathrm{E}_p(|C|_{H_2^j}))
+(\chi_1(p)-\mathrm{E}_p(|C|_{H_1^j}))(\chi_2(p)-1)
\leq\\\leq C''e^{-\gamma'j}
\end{multline*}

Since $p_c(G_1*G_2)\leq \min\{p_c(G_i)\}$ we can take the derivative of their cluster sizes and the following derivative $$\left.\frac{d}{dp}(\chi_1(p)-1)(\chi_2(p)-1)\right|_{p=p_c(G_1*G_2)}
$$
 is positive. There exist $\epsilon$ and $\delta>0$ such that the derivative is bigger than $\delta$ on the interval $[p_c(G_1*G_2),p_c(G_1*G_2)+\epsilon]$.
Let $p_0:=p_c(G_1*G_2)+C''\delta^{-1} e^{-\gamma j}$ then for $j$ large enough $p_0\in[p_c(G_1*G_2),p_c(G_1*G_2)+\epsilon]$ and
\begin{align*}
(\mathrm{E}_{p_c(G_1*G_2)}(|C|_{H_1^j})-1)(\mathrm{E}_{p_c(G_1*G_2)}(|C|_{H_2^j})-1)&\leq
(\mathrm{E}_{p_c(G_1*G_2)}(|C|_{G_1})-1)(\mathrm{E}_{p_c(G_1*G_2)}(|C|_{G_2})-1)= 1\\
(\mathrm{E}_{p_0}(|C|_{H_1^j})-1)(\mathrm{E}_{p_0}(|C|_{H_2^j})-1)&\geq (\mathrm{E}_{p_0}(|C|_{G_1})-1)(\mathrm{E}_{p_0}(|C|_{G_2})-1)- C''e^{-\gamma'j}\\
 &\geq 1+\delta (p_c(G_1*G_2)-p_0)-C''e^{-\gamma'j}
 \geq 1
\end{align*}
Therefore $p_c(H_1^j*H_2^j)\in [p_c(G_1*G_2),p_0]$ and $0\leq p_c(H_1^j*H_2^j)-p_c(G_1*G_2)\leq C''\delta^{-1} e^{-\gamma j}
\stackrel{\scriptscriptstyle j\to\infty}{\longrightarrow}0$.
Thus $p_c(G_1*G_2)=\lim_{j\rightarrow \infty}  p_c(H_1^j*H_2^j)$
which completes the proof.
\end{proof}

\section{Examples}

In Section \ref{recur}, we have presented the expression (\ref{eq:cyc}) of the walk through function $g_i$
for finite cyclic groups. Another relatively easy case is a free group. In
order to find $p_c$ of the free product where one factor is a cyclic
group or a free group, we need to know the expected cluster size.

\begin{proposition}\label{propi}
The expected cluster size in the Cayley graph of the cyclic group
$C_m$ or the free group $F_n$ with respect to a standard set of
generators is given by the following formula
%, respectively
 (in the case of the free group the formula holds for $p<p_c=\frac{1}{2n-1}$)
 \begin{align*}
\mathrm{E}_p(|C|_{C_m})&=\frac{1+p}{1-p}-\frac{p^{m}(m+1)-p^{m+1}(m-1)}{1-p},
\\
\mathrm{E}_p(|C|_{F_n})&=\frac{1+p}{1-(2n-1)p}.
\end{align*}

\end{proposition}

\begin{proof}
The expression for the cluster size in a cyclic group can be
obtained by taking the derivative of $g_i$ given by (\ref{eq:cyc})
and evaluating it at $t=0$.
\begin{align*}
    \mathrm{E}_p(|C|_{C_m})&=1+\sum_{j=1}^{m-1} (j(j-1)(1-p)^2p^{j-1})+  (m(1-p)+p)(m-1)p^{m-1}\\
&=1+\frac{(m^2(p^2-2p+1)-m(3p^2+4p-1)+2p^2)p^{m-1}-2p}{p-1}+
(m(1-p)+p)(m-1)p^{m-1}\\
&=\frac{1+p-m(p^{m}-p^{m+1})-p^m-p^{m+1}}{1-p}
\end{align*}

In order to evaluate the mean cluster size in free group the following observation will be useful. The size of the cluster can be viewed as a sum of indicators that a given vertex is connected to the origin. And the expectation of a sum is a sum of expectations, therefore:
\begin{align}\label{eq:sum}
\mathrm{E}_p(|C|_{G})=\sum_{x\in G}\mathrm{P}_p(x \text{ is
connected to the origin }o).
\end{align}
In the tree, the probability that $x$ is connected to $o$ is equal
to $p^d$, where $d$ is the distance between $x$ and $o$. The
number of vertices at a given distance $d$ is equal to
$2n(2n-1)^{d-1}$. We plug it into the formula (\ref{eq:sum})
and a summation of geometric series completes the
proof.
\begin{align*}
    \mathrm{E}_p(|C|_{F_n})&=1+\sum_{d=1}^{\infty}2n(2n-1)^{d-1} p^{d}=\frac{1+p}{1-(2n-1)p}
\end{align*}

\end{proof}

Using Proposition \ref{propi}, it is easy to compile the following
table of critical probabilities of free products of cyclic groups of
orders $m$ and $n$ (indicated by row and column respectively). For
example the $p_c$ of $PSL(2,\mathbb{Z})=C_2*C_3$ is the unique root of the
 polynomial $2p^2(2-p)(1+p-p^2)-1=2p^5-6p^4+2p^3+4p^2-1$ which lies in the interval $(0,1)$.

\begin{table*}[ht]
\centering
        \begin{tabular}{c|c|c|c|c|c|c|c|}
            $m\setminus n$&2&3&4&5&10&100&$\infty$\\
            \hline
            2&1&.5199&.4613&.4414&.4271&.4268&.4268\\
            \hline
            4&.4613&.3754&.3539&.3468&.3427&.3426&.3426\\
            \hline
            10&.4271&.3605&.3427&.3367&.3334&.3334&.3334\\
            \hline
        \end{tabular}
  \caption{Values of $p_c$ for some $C_m*C_n$}
  \label{tab}
\end{table*}

We can see that the critical probability of the Cayley graph of
$C_m*C_n$ decreases as the orders grow and the limit of $p_c(C_m*C_n)$ as $m,n\to \infty$
is $p_c(F_2)=\frac13$. It follows from the fact
that $\mathrm{E}_p(|C|_{C_m})$ grows with the order of the cyclic
group and the limit is equal to the one for the infinite cyclic
group, in particular for infinite cyclic  group we have
$\mathrm{E}_p(|C|_\mathbb{Z})=(1-p)^2\sum_{i=1}^{\infty}i^2p^{j-1}=(1+p)/(1-p)$.
A more general situation is treated by the Theorem \ref{limit}.

It is well known and follows from formula (\ref{propi}) that  the critical probability for a regular tree 
(the Cayley graph of the free group of rank $n$) is equal to $\frac{1}{2n-1}$. The inequality  for $p_c$
involving the Cheeger constant and inequality (\ref{eq:est2}) become equalities in this
case since $h(F_n)=2n-2$. These two inequalities become strict even if we consider free products of finite cyclic groups (which are virtually free).

Note that if a group is virtually $\Z$ then $p_c=1$ and $h(\G)=0$, so  we again have equality $p_c=\frac{1}{h(\G)+1}$ (the other inequality  is strict if the generating set contains more than one element). An interesting problem is whether  the inequality $p_c\le \frac{1}{h(\G)+1}$ is strict except for free groups with a free generating set. The next proposition shows that even for virtually free groups the inequality can be strict.

\begin{proposition}\label{cheeger}
Consider the free product of two finite cyclic groups $C_m*C_n$,
such that $(m-1)(n-1)>1$ with the natural generators. Then the inequalities \eqref{eq:est1},
\eqref{eq:est2} are strict, that is
\begin{align*}
\frac{1}{3}<p_c<\frac{1}{h(C_m*C_n)+1}.
\end{align*}
\end{proposition}

\begin{proof}
The first inequality follows from the observation that
$\mathrm{E}_p(|C|_\mathbb{Z})>\mathrm{E}_p(|C|_{C_m})$ for all $m$
and $p>0$. For the second inequality we will show that
\begin{align}
h(C_n*C_m)\leq
2-\max\left(\frac{2m}{n(m-1)},\frac{2n}{m(n-1)}\right)<1/p_c-1.\label{eq:cheeg}
\end{align}
Consider a set $S_1$ consisting of a cycle of $C_n$ at the origin in the Cayley graph of $C_n*C_m$. We construct $S_k$
inductively by including the whole cycle $C_m$ at every point of the
boundary of $S_{k-1}$ and then adding the next generation of cycles $C_n$ to the
$S_{k-1}$.

Then
\begin{align*}
|\partial S_1|&=2n\\
|S_1|&=n\\
 |\partial S_{k+1}|&=(n-1)(m-1)|\partial S_k|\\
 &=2n((n-1)(m-1))^k
 \\
 |S_{k+1}|&=|S_k|+|\partial S_{k}|(m-1)n\\
 &=n+n^2(m-1)\sum_{i=0}^{k-1}((n-1)(m-1))^i\\
 &=n+n^2(m-1)\frac{((n-1)(m-1))^k-1}{(n-1)(m-1)-1}\\
\frac{|\partial S_{k+1}|}{|S_{k+1}|} &=
\frac{2n((n-1)(m-1))^k(mn-m-n)}{n(mn-m-n)+n^2(m-1)((n-1)(m-1))^k-1)}\\
&\stackrel{k\to \infty}{\longrightarrow}
  \frac{2(mn-m-n)}{n(m-1)}\\
h(C_n*C_m) &\leq 2-\frac{2m}{n(m-1)}
\end{align*}

Assume that $n\leq m$ and set
\begin{align}\label{eq:odh}
p_1=\frac{1}{3-\frac{2}{n}}\leq \frac{1}{3-\frac{2m}{n(m-1)}}\leq
\frac{1}{h(C_n*C_m)+1}.
\end{align}

It is enough to prove $p_c<p_1$. By Theorem \ref{fortwo}, it means
we need to show that
$$(\mathrm{E}_{p_1}(|C|_{C_m})-1)(\mathrm{E}_{p_1}(|C|_{C_n})-1)>1.$$

We have:
\begin{align*}
    (\mathrm{E}_{p_1}(|C|_{C_n})-1)
&   =
\frac{2p_1+p_1^{n+1}(n-1)-p_1^{n}(n+1)}{1-p_1}\\
&=1+\frac{1-(n^2+n-1)\left(\frac{n}{3n-2}\right)^n}{n-1}
\end{align*}

Now for $n\geq 3$:
\begin{align*}
(n^2+n-1)\left(\frac{1}{3}\left(1+\frac{2}{3n-2}\right)\right)^n
\leq(n^2+n-1)\left(\frac{3}{7}\right)^n<1.
\end{align*}
Thus if $ 3\leq n\leq m$ we have
\begin{align*}
    (\mathrm{E}_{p_1}(|C|_{C_m})-1)(\mathrm{E}_{p_1}(|C|_{C_n})-1)\geq
    (\mathrm{E}_{\frac{m}{3m-2}}(|C|_{C_m})-1)(\mathrm{E}_{\frac{n}{3n-2}}(|C|_{C_n})-1)
        >1.
\end{align*}

It remains to consider the case $n=2$. If $m=3$ we will use the
stronger estimate in (\ref{eq:odh})
$p_2:=\frac{1}{3-\frac{2m}{n(m-1)}}=\frac{3}{2}$; otherwise we use
again $p_1$, now equal to $\frac{1}{2}$. By plugging in we obtain
\begin{align*}
        (\mathrm{E}_{p_2}(|C|_{C_2})-1)(\mathrm{E}_{p_2}(|C|_{C_3})-1)&=1.4\\
      (\mathrm{E}_{p_1}(|C|_{C_2})-1)(\mathrm{E}_{p_1}(|C|_{C_m})-1)&\geq
      (\mathrm{E}_{p_1}(|C|_{C_2})-1)(\mathrm{E}_{p_1}(|C|_{C_4})-1)=1.2
\end{align*}
Therefore in all cases we obtained that $p_c<1/(h(C_m*C_n)+1)$, as
required.
\end{proof}

\comment
\section{The critical exponent}

In this section we are going to study the critical exponents of
free products, in particular we will prove that $\beta$ is equal to $1$
 under conditions specified in Theorem \ref{expon}.

Recall that $\beta$ and $\gamma$ are called critical
exponents if the limits
\begin{align*}
\lim_{p\to p_{c+}}\frac{\theta(p)}{(p-p_c)^\beta}\text{, and }
\lim_{p\to p_{c-}}\frac{\mathrm{E}_{p_c}(|C|_{G_1*G_2})}{(p_c-p)^\gamma}
\end{align*}
exist and are nonzero, respectively.

\begin{proposition}\label{exptwo}
Let $G_1=\langle S_1\rangle$,
$G_2=\langle S_2\rangle$ be two non trivial finitely generated groups such that $|G_1|>2$ or $|G_2|>2$. Consider
the Cayley graph of the free product $G_1 * G_2$ with respect to the
generating set $S_1\cup S_2$. Then the critical exponent $\gamma=-1$.
If moreover $\mathrm{E}_{p_c}(|C|_{G_i}^2)<\infty$ at $p_c=p_c(G_1*G_2)$ for $i=1,2$ then the critical exponent $\beta=1$.
\end{proposition}

\begin{proof}
In order to prove that $\gamma=-1$ we will show that the following derivative from the left exists and is finite. Denote the left derivative (in $p$) by $\frac{d}{d_{-}p}$.
\begin{align*}
\lim_{p\to p_{c-}}\frac{(\mathrm{E}_{p_c}(|C|_{G_1*G_2}))^{-1}}{p_c-p}=-\left.\frac{d}{d _-p}(\mathrm{E}_{p_c}(|C|_{G_1*G_2}))^{-1}\right|_{p=p_c}
\end{align*}
Since $p_c(G_1*G_2)<p_c(G_i)$, $\mathrm{E}_{p}(|C|_{G_i})$ is finite
by Lemma \ref{lemma} and has finite nonzero derivative at
$p=p_c(G_1*G_2)$. Therefore
\begin{multline*}
\left.\frac{d}{d _-p}(\mathrm{E}_{p_c}(|C|_{G_1*G_2}))^{-1}\right|_{p=p_c}
=\left.\frac{d}{d _-p}
\frac{\mathrm{E}_{p}(|C|_{G_1})+\mathrm{E}_{p}(|C|_{G_2})-\mathrm{E}_{p}(|C|_{G_1})\mathrm{E}_{p}(|C|_{G_2})}{\mathrm{E}_{p}(|C|_{G_1})\mathrm{E}_{p}(|C|_{G_2})}\right|_{p=p_c}\\
=\frac{\left.\frac{d}{d _-p} \mathrm{E}_{p}(|C|_{G_1}) \right|_{p=p_c}(1-\mathrm{E}_{p_c}(|C|_{G_2}))+\left.\frac{d}{d _-p} \mathrm{E}_{p}(|C|_{G_2}) \right|_{p=p_c}(1-\mathrm{E}_{p_c}(|C|_{G_1}))   }{\mathrm{E}_{p_c}(|C|_{G_1})\mathrm{E}_{p_c}(|C|_{G_2})}
<0
\end{multline*}
Thus the required limit for $\gamma=-1$ exists and is finite nonzero.

Next, for proving that $\beta=1$ we need to show that the right
derivative of $\theta$ at $p_c$ is finite and nonzero.
First observe that the walk through functions $g_i$ are power series (in $X$) with a zero constant term. So we can write:
\begin{align*}
g_i(p,t)&=\sum_{j=1}^{|G_i|-1}\frac{1}{j!}\left.\frac{\partial^j g_i(p,t)}{\partial t^j}\right|_{t=0}t^j,\\
\intertext{where}\qquad
\left.\frac{\partial^k g_i(p,t)}{\partial t^k}\right|_{t=0}&=\mathrm{E}_p[(|C|_{G_i}-1)\dots(|C|_{G_i}-k)].
\end{align*}
Now we plug these expressions into the formula
$B=g_2(p,g_1(p,B))$. Since $p_c(G_1*G_2)<p_c(G_i)$ we can consider $p_c(G_1*G_2)<p<\min\{p_c(G_i)\}$ for which $\chi_i(p)$ is analytical. Moreover $B(p)>0$ and we can
cancel $B$ in the equality as follows:
\begin{align*}
B(p)&=\sum_{j=1}^{|G_2|-1}\frac{1}{j!}\left.\frac{\partial^j g_2(p,t)}{\partial t^j}\right|_{t=0}\left(\sum_{i=1}^{|G_1|-1}\frac{1}{i!}\left.\frac{\partial^i g_1(p,t)}{\partial t^i}\right|_{t=0}B(p)^i\right)^j\\
1&=\left.\frac{\partial g_2(p,t)}{\partial t}\right|_{t=0}\left.\frac{\partial g_1(p,t)}{\partial t}\right|_{t=0}+
\frac{B(p)}{2}\left(\left.\frac{\partial g_2(p,t)}{\partial t}\right|_{t=0}\left.\frac{\partial^2 g_1(p,t)}{\partial t^2}\right|_{t=0}\right.
\\&
\quad +\left.\left.\frac{\partial^2 g_2(p,t)}{\partial t^2}\right|_{t=0}\left(\left.\frac{\partial g_1(p,t)}{\partial t}\right|_{t=0}\right)^2\right)+B(p)^2(O(1))\\
\frac{1-\left.\frac{\partial g_2(p,t)}{\partial t}\right|_{t=0}\left.\frac{\partial g_1(p,t)}{\partial t}\right|_{t=0}}{p-p_c}&=\frac{B(p)}{2(p-p_c)}\left(\left.\frac{\partial g_2(p,t)}{\partial t}\right|_{t=0}\left.\frac{\partial^2 g_1(p,t)}{\partial t^2}\right|_{t=0}\right.
\\&\quad\left.
+\left.\frac{\partial^2 g_2(p,t)}{\partial t^2}\right|_{t=0}\left(\left.\frac{\partial g_1(p,t)}{\partial t}\right|_{t=0}\right)^2\right)
+\frac{B(p)^2}{p-p_c}(O(1))\\
\end{align*}

As $p$ approaches $p_c$, $B(p)$ goes to $0$ and so does the last
term in the above equation. Denote $E_p(|C|_{G_i})$ by $\chi_i(p)$. This leaves us with the
following expression of the derivative of $B$ from the right at
$p_c$:
\begin{align*}
\left.\frac{d B(p)}{d_{+} p}\right|_{p=p_{c}}&=
\frac{2
\left.\frac{d}{d_+ p}\left(1-\left.\frac{\partial g_2(p,t)}{\partial t}\right|_{t=0}\left.\frac{\partial g_1(p,t)}{\partial t}\right|_{t=0}\right)\right|_{p=p_c}
}{
\left.\frac{\partial g_2(p_c,t)}{\partial t}\right|_{t=0}\left.\frac{\partial^2 g_1(p_c,t)}{\partial t^2}\right|_{t=0}+\left.\frac{\partial^2 g_2(p_c,t)}{\partial t^2}\right|_{t=0}\left(\left.\frac{\partial g_1(p_c,t)}{\partial t}\right|_{t=0}\right)^2}\\
&=\frac{-2
\left.\frac{d}{d_+ p}\left((\chi_1(p)-1)(\chi_2(p)-1)\right)\right|_{p=p_c}
}{\mathrm{E}_{p_c}[(|C|_{G_1}-1)(|C|_{G_1}-2)](\chi_2(p_c)-1)+\mathrm{E}_{p_c}[(|C|_{G_2}-1)(|C|_{G_2}-2)](\chi_1(p_c)-1)^2}
\end{align*}
As soon as  $(|G_1|-1)(|G_2|-1)>1$ the denominator is nonzero. And if $\mathrm{E}_{p_c}(|C|_{G_i}^2)<\infty$ the denominator is finite,  the
numerator equals to the derivative (in $p$) of
$1-(\chi_1(p)-1)(\chi_2(p)-1)$ which is finite and nonzero.  Thus the derivative from the right of $B(p)$ at $p_c$ is nonzero and finite.

Recall the expression for percolation function:
$\theta(p)=B(p)+g_1(p,B(p))-B(p)g_1(p,B(p))$ (formula (\ref{eq:theta})).
Note that $B(p_c)=0$. Thus for the derivative of $\theta$ we have:
\begin{align}
\left.\frac{d\theta(p)}{d_+ p}\right|_{p=p_c}
&=
\left.\frac{d\left(B(p)+g_1(p,B(p))(1-B(p)) \right)}{d_+ p}\right|_{p=p_c}\notag\\
&=
\left.\frac{d B(p)}{d_+ p}\right|_{p=p_c}-\left. g_1(p,0)\frac{d B(p)}{d_+ p}\right|_{p=p_c}
+\left.\frac{\partial g_1(p,t)}{\partial t}\right|_{t=B(p)}\left.\frac{d B(p)}{d_+ p}\right|_{p=p_c}\left(1-B(p_c)\right)
\notag
\\&
\quad+\left.\frac{\partial g_1(p,B(p_c))}{\partial p}\right|_{p=p_c}(1-B(p_c))
\notag
\\
&=
\left.\frac{d B(p)}{d_+ p}\right|_{p=p_c}\left( 1+ \left.\frac{\partial g_1(p_c,t)}{\partial t}\right|_{t=0} \right)
\notag
\\
\left.\frac{d\theta(p)}{d_+ p}\right|_{p=p_c}
&= \, \delta \left.\frac{d}{d_+ p}\left((\chi_1(p)-1)(\chi_2(p)-1)\right)\right|_{p=p_c},\label{eq:dth}
\\
\intertext{where}
\qquad
\delta
&=\frac{-2\chi_1(p_c)}
{(\mathrm{E}_{p_c}(|C|_{G_1})^2-1)(\chi_2(p_c)-1)+
(\mathrm{E}_{p_c}(|C|_{G_2})^2-1)(\chi_1(p_c)-1)^2-
3\chi_1(p_c)
}\label{eq:por1}\\
&=\frac{-2\chi_2(p_c)}
{(\mathrm{E}_{p_c}(|C|_{G_1})^2-1)(\chi_2(p_c)-1)^2+
(\mathrm{E}_{p_c}(|C|_{G_2})^2-1)(\chi_1(p_c)-1)-
3\chi_2(p_c)}
\label{eq:por2}
\end{align}
%&=-2\mathrm{E}_{p_c}[|C|_{G_2}]/ (
%\mathrm{E}_{p_c}[(|C|_{G_2})^2-1]\mathrm{E}_{p_c}[(|C|_{G_1})-1]
%\notag
%\\&
%\quad+\mathrm{E}_{p_c}[(|C|_{G_1})^2-1](\mathrm{E}_{p_c}[(|C|_{G_2})-1])^2
%-3\mathrm{E}_{p_c}[|C|_{G_1}])
%\label{eq:por1}
%\\&=
%-2\mathrm{E}_{p_c}[|C|_{G_1}]/ (
%\mathrm{E}_{p_c}[(|C|_{G_1})^2-1]\mathrm{E}_{p_c}[(|C|_{G_2})-1]
%\notag
%\\&
%\quad+\mathrm{E}_{p_c}[(|C|_{G_2})^2-1](\mathrm{E}_{p_c}[(|C|_{G_1})-1])^2
%-3\mathrm{E}_{p_c}[|C|_{G_2}])\label{eq:por2}
The equalities (\ref{eq:por1}) and (\ref{eq:por2}) correspond to
exchanging $G_1$ with $G_2$ and can be obtained from each other
using the fact that
$(\chi_1(p_c)-1)(\chi_2(p_c)-1)=1$.
Again the derivative of
$\theta$ at $p_c$ from the right exists and is nonzero and finite
(because the derivative of $B$ at $p_c$ is such). Therefore the
critical exponent is equal to $1$.\end{proof}

\begin{proof}[Proof of Proposition \ref{expon}]
The results in the Proposition \ref{exptwo} can be again inductively generalized to the free product of $n$ (finite) groups.
If $G_i$ is finite then clearly $\mathrm{E}_{p_c}(|C|_{G_i}^2)<\infty$ and the only thing we need to show is that this second moment remains finite under taking free products.

We can express $\mathrm{E}_p(|C|^2_{G_1*G_2})$ in terms of $\chi_i(p)$ and $\mathrm{E}_p(|C|^2_{G_i})$. Recall the splitting of the Cayley graph (without origin) into two parts: part $P_1$ consists of vertices $v$ such that any simple path from the
origin to $v$ first visits a vertex corresponding to an element of
$G_1$ (in fact $S_1$), part $P_2$ is defined similarly for
$G_2$. Recall the notation $|C|_{P_i}=|C\cap P_i|_{G_1*G_2}$. Then
\begin{align*}
\mathrm{E}_p(|C|^2_{G_1*G_2})&=2(\mathrm{E}_p(|C|_{P_1})+1)(\mathrm{E}_p(|C|_{P_2})+1)-1+\mathrm{E}_p(|C|^2_{P_1})+\mathrm{E}_p(|C|^2_{P_2})\\
\mathrm{E}_p(|C|^2_{P_1})&=(\chi_1(p)-1)(\mathrm{E}_p(|C|^2_{P_2})+2\chi_2(p)-1)\\
&\quad +(\mathrm{E}_p(|C|^2_{G_1})-2\chi_1(p)+1)(\mathrm{E}_p(|C|_{P_2})+1)^2
\end{align*}

A similar formula holds for $\mathrm{E}_p(|C|^2_{P_2})$. We can compose them into one equation and solve it for $\mathrm{E}_p(|C|^2_{P_i})$. Recall that in Section \ref{ctyri} we have found $\mathrm{E}_p(|C|_{P_i})$.

The resulting expression for $\mathrm{E}_p(|C|^2_{G_1*G_2})$ is then
\begin{align*}
\mathrm{E}_p(|C|^2_{G_1*G_2})=\frac{\gamma}{(\chi_1(p)+\chi_2(p)-\chi_1(p)\chi_2(p))^3}
\end{align*}
where
\begin{align*}
\gamma&=
6\chi_2(p)\chi_1^2(p)
+6\chi_2^2(p)\chi_1(p)
+10\chi_2^2(p)\chi_1^3(p)%\\&\quad
+10\chi_2^3(p)\chi_1^2(p)
-16\chi_2^2(p)\chi_1^2(p)
-7\chi_2(p)\chi_1^3(p)\\&\quad
-7\chi_2^3(p)\chi_1(p)
-5\chi_2^3(p)\chi_1^3(p)
+\chi_1^3(p) +\chi_2^3(p)%\\&\quad
+\chi_2^3(p)\chi_1(p^2)
+\chi_1^3(p)\chi_2(p^2).
\end{align*}
Therefore $\mathrm{E}_p(|C|^2_{G_1*G_2})$ is finite for all $p<p_c$. This completes the proof.
\end{proof}

\endcomment

%\bibliography{biblio}
%\bibliographystyle{amsplain}

\providecommand{\bysame}{\leavevmode\hbox to3em{\hrulefill}\thinspace}
\providecommand{\MR}{\relax\ifhmode\unskip\space\fi MR }
\providecommand{\MRhref}[2]{%  
\href{http://www.ams.org/mathscinet-getitem?mr=#1}{#2}}
\providecommand{\href}[2]{#2}

\noindent Iva Koz\'akov\'a\\
Department of Mathematics\\
Vanderbilt University\\
iva.kozakova@vanderbilt.edu

\end{document}